\documentclass[11pt,fleqn]{article}
\usepackage{geometry}
\usepackage{amssymb,amsmath}
\usepackage{amsthm}
\usepackage[square,sort,comma,numbers]{natbib}
\usepackage{graphicx, float}
\usepackage[final]{showkeys}
\usepackage{breqn}
\theoremstyle{plain}

\theoremstyle{definition}
\newtheorem{exmp}{Example} 

\makeatletter
\newdimen\mymathindent
\newenvironment{bulletequation}%
{\@beginparpenalty\predisplaypenalty
	\@endparpenalty\postdisplaypenalty
	\refstepcounter{equation}%
	\trivlist \item[]\leavevmode
	\hb@xt@\linewidth\bgroup $\m@th
	\displaystyle
	\hskip\mymathindent}%
{$\hfil 
	\displaywidth\linewidth\hbox{\@eqnnum}%
	\egroup
	\endtrivlist}
\makeatother
\begin{document}
\newcommand*{\Scale}[2][4]{\scalebox{#1}{\ensuremath{#2}}}
	\newcommand{\inorm}[1]{\lVert #1\rVert_{\infty}}
	\newcommand{\mnorm}[1]{\lvert #1\rvert}
	\newcommand{\ml}[2]{E_{#1}(t^{#1} #2)}
	\newcommand{\gml}[3]{E_{#1}(t^{#2} #3)}
	\newcommand{\Gml}[3]{E_{#1}({#2}\:#3)}
	\newcommand{\pder}[2]{\dfrac{\partial^{#1}}{\partial #2^{#1}}}
	\newcommand{\norm}[1]{\lVert #1\rVert}
	\newcommand{\mll}[2]{E_{#1}(#2 t^{#1})}
\begin{center}
	{\LARGE Local Stable Manifold theorem for fractional systems revisited}
\end{center}
\vspace{0.5cm}
\begin{center}
	{
		{\Large Amey Deshpande\\}
		\textit{ Department of Mathematics, Savitribai Phule Pune University, Pune - 411007, India,\\ 2009asdeshpande@gmail.com\\}
	}
	{\Large Varsha Daftardar-Gejji\\}
	\textit{Department of Mathematics, Savitribai Phule Pune University, Pune - 411007, India\\ vsgejji@math.unipune.ac.in, vsgejji@gmail.com}\\
	\end{center}
	\begin{abstract}
	The subject of fractional calculus has witnessed rapid development over past few decades. In particular the area of fractional differential equations has received considerable attention. Several theoretical results have been obtained and powerful numerical methods have been developed. In spite of the extensive numerical simulations that have been carried out in the area of fractional order dynamical systems, analytical results obtained are very few. In pursuance to this, present authors have extended local stable manifold theorem in case of fractional systems \cite{deshpande2016local}. 	Cong et al. \cite{cong2016stable} have pointed out discrepancies in the asymptotic expansion of two-parameter Mittag-Leffler functions with matrix argument (\textit{cf.} Lemma 4 part 2 of article \cite{deshpande2016local}). In the present communication we give the corrected expansion of the same and prove the local stable manifold theorem by following the same approach given in \cite{deshpande2016local}.
	\end{abstract}
	\section{Introduction}
	The area of fractional order dynamical systems (FODS) is being actively pursued in the last few years due to their applications in diverse fields including viscoelasticity \cite{mainardi2010fractional}, fractional order control systems \cite{podlubny1999fractional}, mechanics \cite{riewe1996nonconservative, riewe1997mechanics}, bioengineering \cite{magin2004fractional}, economics \cite{jun2001study} and so on. The study of stability of FODS was initiated by Matignon in the year 1996 \cite{matignon1996stability}. Further important development in this area was due to Grigorenko and Grigorenko \cite{grigorenko2003chaotic} who studied fractional Lorenz system and proved that the order of the fractional derivatives act as a chaos controlling parameter. Since then a lot of simulation work has been carried out to explore various FODS \cite{tavazoei2008chaotic, daftardar2010chaos, bhalekar2010fractional, li2006chaos, yu2008synchronization}. In spite of  extensive numerical work, analytical results obtained in this area are very few. Concept of flow is pivotal in the theory of dynamical systems. One of the major hurdle in the development of FODS is that their solutions do not satisfy semi-group property. As a consequence, concept of `flow' in fractional framework has not been satisfactorily formulated. 
	
	The local stable manifold theorem is one of the basic results in the realm of the dynamical systems. Present authors have generalized the local stable manifold for fractional systems. A proof of this theorem was presented in \cite{deshpande2016local}.
	

	Cong et al. \cite{cong2016stable} have pointed out discrepancies in asymptotic expansion of two-parameter Mittag-Leffler function with matrix argument(Lemma 4, part 2) of the proof of local stable manifold theorem given by the present authors in \cite{deshpande2016local}. In the present paper we provide correct formulation of the Lemma 4 part 2, and consequent changes in some expressions of Lemma 5, 6, 8 and the Step II of the proof of the main theorem given in \cite{deshpande2016local}. Further we include the discussion about examples discussed both in \cite{deshpande2016local} and \cite{cong2016stable}. Thus the proof of the local stable manifold theorem given by us in \cite{deshpande2016local} continues to hold true in view of these corrections.
	\section{Corrections}
	The corrected version of Lemma 4 Part 2 of the ref. \cite{deshpande2016local} is presented below. Note that the Lemma 4 part 1 remains as it is.\\
	\textbf{Lemma 4 part 2.}
		For $0 < p < 1$, $j=1,2 \cdots, l$ and $q \in \mathbb{N}\backslash\{1\}$.
		\begin{equation}
			\gml{p,p}{p}{J_j} = t^{-p} \widetilde{B_j(t)} + \widetilde{C_j(t)},
		\end{equation} where $\widetilde{B_j(t)}$ and $\widetilde{C_j(t)}$ are $n_j \times n_j$ matrices defined as 
		\begin{align} \label{mat6}
			\nonumber \widetilde{B_j(t)} &:= \hfill 0    \hfill & j=1,2,\cdots,s, \\
			\widetilde{B_j(t)}&:= \hfill \begin{pmatrix}
				\widetilde{\Psi}_0(t,\lambda_j)& \widetilde{\Psi}_1(t,\lambda_j) & \cdots & \widetilde{\Psi}_{n_j-1}(t,\lambda_j)\\
				& \widetilde{\Psi}_0(t,\lambda_j) & \widetilde{\Psi}_1(t,\lambda_j)  \cdots & \\
				& & & \vdots \\
				&  & \ddots  & \\
				&  &  & \widetilde{\Psi}_0(t,\lambda_j)
			\end{pmatrix}, \hfill & j=s+1,s+2,\cdots,l, 
		\end{align}
		
		\begin{align} \label{mat4}
			\widetilde{C_j(t)} := \begin{pmatrix}
				\widetilde{\Delta}_0(t,\lambda_j) & \widetilde{\Delta}_1(t,\lambda_j) & \cdots & \widetilde{\Delta}_{n_j-1}(t,\lambda_j)\\
				&  \widetilde{\Delta}_0(t,\lambda_j)& \widetilde{\Delta}_1(t,\lambda_j)  \cdots & \\
				& & & \vdots \\
				&  & \ddots  & \\
				&  &  & \widetilde{\Delta}_0(t,\lambda_j)
			\end{pmatrix}_{n_j \times n_j}, ~~ j=1,2,\cdots,l,
		\end{align} 
		and 
		\begin{equation}
		\widetilde{\psi_m} (t, \lambda_j) = \frac{1}{m!} \left[ \pder{m+1}{\lambda} \exp(t \lambda_j^{\frac{1}{p}})\right],~~m=0,1,2,\cdots,n_j-1. 
		\end{equation}
\begin{align} \label{delta}
\widetilde{\Delta_m}(t,\lambda_j):= \dfrac{1}{m!}\bigg( \sum_{k=2}^{q} \dfrac{(-1)^{m+2}(k+m)!}{(k-1)!} \dfrac{ \lambda_j ^{-k-m-1} t^{-pk}}{\Gamma(1-pk)} + O( \lvert\lambda_{j}\rvert^{-q-m} t^{-p-pq}) \bigg). 
\end{align}
Let $\widetilde{B(t)}$ and $\widetilde{C(t)}$ denote the block diagonal matrices consisting of $\widetilde{B_j(t)}$ and $\widetilde{C_j(t)}$ on the diagonal respectively. 
		
Then $\gml{p,p}{p}{A} = t^{-p} \widetilde{B(t)} + \widetilde{C(t)}$, where $J_j$ and $A$ are as defined in the paper \cite{deshpande2016local}.

	\begin{proof}
		From \citep[thm 1.3]{podlubny1999fractional} for sufficiently large $t$,
		\begin{align}
		\nonumber\gml{p,p}{p}{\lambda_j} &= \frac{1}{p} t^{1-p} \lambda_j ^{\frac{1}{p}-1} \exp(t \lambda_j ^{\frac{1}{p}}) - \sum_{k=1}^{q} \frac{t^{-pk} \lambda_j ^{-k}}{\Gamma(p-pk)} + O(\mnorm{\lambda_j} ^{-q-1}) \\
		\nonumber&= t^{-p} \left( \frac{1}{p} t \lambda_j ^{\frac{1}{p}-1} \exp(t \lambda_j ^{\frac{1}{p}})\right) - \sum_{k=1}^{q} \frac{t^{-pk} \lambda_j ^{-k}}{\Gamma(p-pk)} + O(\mnorm{\lambda_j}) \\
		& = t^{-p} \pder{ }{\lambda} \exp(t \lambda_j ^{\frac{1}{p}}) - \sum_{k=1}^{q} \frac{t^{-pk} \lambda_j ^{-k}}{\Gamma(p-pk)} + O(\mnorm{\lambda_j}). \label{err1} 
		\end{align} 
		Differentiating eqn. \eqref{err1} $m$-times with respect to $\lambda$ and multiplying by $\frac{1}{m!}$ we obtain
		\begin{align}
		\nonumber\frac{1}{m!} \pder{m}{\lambda} \gml{p,p}{p}{\lambda_j} & = \frac{t^{-p}}{m!} \pder{m+1}{\lambda} \exp(t \lambda_j ^{\frac{1}{p}}) + \widetilde{\Delta_m}(t,\lambda_j) \\
		&= t^{-p} \widetilde{\Psi_m}(t, \lambda_j) + \widetilde{\Delta_m}(t,\lambda_j), ~~m=0,1,\cdots,n_j-1. \label{err2}  
		\end{align}
		The required results follow from the eqn. \eqref{err2}.
	\end{proof}
	
	\textbf{The asymptotic behavior of $\widetilde{C(t)}$ and $\widetilde{B(t)}$:} Along the same lines as in \citep[Lemma 5]{deshpande2016local}, we get
	\begin{equation}
    \mnorm{\widetilde{B(-t)}} \leq \sum_{m=0}^{n_{\tilde{r}}-1} \frac{1}{p^{m+1}} M(\lambda_{\tilde{r}},m+1) t^{m+1}  e^{-t \alpha}.  
	\end{equation}
	
	Further it is noted that $\mnorm{C(t)}$ and $\mnorm{\widetilde{C(t)}}$ have the same asymptotic behavior. Since 
	\begin{equation}
	\mnorm{\widetilde{C(t)}} \leq \tilde{K_1}(q, \lambda_{\tilde{r}}) t^{-2p} + \tilde{K_2}(q, \lambda_{\tilde{r}}) t^{-p-pq},
	\end{equation}
	where $\tilde{K_1} (q, \lambda_{\tilde{r}})$ and $\tilde{K_2}(q, \lambda_{\tilde{r}})$ denote the arbitrary constants depending on $q$ and $\lambda_{\tilde{r}}$.
	
	
	The revised version of \citep[Lemma 6 (eqn. (35))]{deshpande2016local} takes the following form:
	
	\textbf{Lemma 6.}\label{lem6}
		For any $t, \tau > 0$,
		\begin{equation}
		(t-\tau) B(t) \widetilde{B(-\tau)} = \frac{- \tau}{p} \widetilde{B(t-\tau)}.
		\end{equation}
		
		\begin{proof}
			Since $B, \widetilde{B}$ are block diagonal matrices, it is sufficient to prove 
			\begin{equation}
			(t-\tau) B_j(t) \widetilde{B_j(-\tau)} = \frac{-\tau}{p} \widetilde{B_j (t-\tau)}, ~~j=1,2,\cdots,l.
			\end{equation}
			
			In view of the definitions of $B_j$ and $\widetilde{B_j}$ it suffices to find the elements in the first row of the matrix $B_j(t) \widetilde{B_j (-\tau)}$. Denote the elements in the first row by $(\eta_0, \eta_1, \eta_2, \cdots, \eta_{n_j -1})$. For $m=1,2, \cdots, n_j -1$, 
			\begin{align}
				(t-\tau) \eta_m &= (t-\tau) \sum_{i=0}^{m} \Psi_i(t) \widetilde{\Psi_{m-i}} (-\tau) \\ &= (t-\tau) \sum_{i=0}^{m} \frac{1}{i!} \pder{i}{\lambda}\left( \frac{1}{p}\exp(t\lambda_j^{\frac{1}{p}})\right) \frac{1}{(m-i)!} \pder{m-i+1}{\lambda}\left( \exp(-\tau\lambda_j^{\frac{1}{p}})\right)\\
				&= \frac{(t-\tau)}{p ~m!}\sum_{i=0}^{m} \binom{m}{i} \pder{i}{\lambda}\left( \exp(t\lambda_j^{1/p})\right) \pder{m-i}{\lambda}\left( \pder{ }{\lambda}\exp(-\tau\lambda_j^{1/p})\right) \\
				&= \frac{(t-\tau)}{p~m!} \pder{m}{\lambda}\left( \exp(t\lambda_j^{1/p})~ \pder{ }{\lambda}\exp(-\tau\lambda_j^{1/p}) \right)\\
				&= \frac{-\tau}{p~m!} \pder{m}{\lambda} \left( \frac{(t-\tau)}{p} \lambda_j ^{\frac{1}{p}-1} \exp((t-\tau) \lambda_j ^{\frac{1}{p}})	\right)\\
				& = \frac{-\tau}{p}\left[ \frac{1}{m!} \pder{m+1}{\lambda} \exp((t-\tau) \lambda_j ^{\frac{1}{p}}) \right] \\
				&= \frac{-\tau}{p} \widetilde{\Psi_m (t-\tau)}. 			
			\end{align} 
		\end{proof}

	In the light of the above corrected Lemma, the operator $T_{\sigma}(x)$ \citep[ eqn.(44) and eqn.(45)]{deshpande2016local} should be replaced by the following equations. Note that \citep[eqn.(45)]{deshpande2016local} remains same. Keeping all notations same as in \cite{deshpande2016local} we define
\begin{align} \label{define}
\begin{split}
\pi_u(T_{\sigma}x(t)) &= \ml{p}{A} ~ \int_{0}^{\infty} \frac{p}{\tau} \widetilde{B(-\tau)}~ \pi_u f(x(\tau))d\tau ~+\\ &\int_{0}^{t} (t-\tau)^{p-1} E_{p,p}( (t-\tau)^p A)~ \pi_u f(x(\tau))d\tau,
\end{split}\\
\pi_s(T_{\sigma}x(t)) &=  E_p(t^p A)\sigma + \int_{0}^{t} (t-\tau)^{p-1} E_{p,p}((t-\tau)^p A) ~\pi_s f(x(\tau))d\tau.
\end{align}

Let $N > 0$, $\tilde{r}>0$ be such that $\mnorm{\gml{p,\beta}{p}{\lambda_{\tilde{r}}}} = \max\limits_{1 \leq j \leq l} \mnorm{\gml{p,\beta}{p}{\lambda_j}},~\beta>0$. Denote $\mnorm{A}$ by $a$ and $q \in \mathbb{N}\backslash\{1\}$. As a consequence of the above corrections, \citep[Lemma 8 part 1 and part 2 (eqn. (46) and eqn. (47))]{deshpande2016local} are revised below while \citep[eqn. (48) and eqn. (49)]{deshpande2016local} remain unaltered.

\textbf{Lemma 8 part 1 and part 2.} \label{lem8}
Let for $g \in C[I,\mathbb{R}^n]$, and $t>N$,
\begin{enumerate}
	\item 
	\begin{bulletequation}
		\left\| \int_{0}^{\infty} \frac{p}{\tau} C(t)\widetilde{B(-\tau)} \; \pi_ug(\tau) d\tau \right\| \leq K_3(N,q,\lambda_{\tilde{r}}) \inorm{g},
	\end{bulletequation}
	where $K_3$ is an arbitrary constant. 
	\item \begin{align}
	\begin{split}
	\left\| \int_{0}^{t} (t-\tau)^{p-1} \Gml{p,p}{(t-\tau)^p}{A}\; \pi_ug(\tau) d\tau + \int_{0}^{\infty}\frac{p}{\tau} B(t)\widetilde{B(-\tau)}~ \pi_ug(\tau) d\tau\right\|\\ < K_5(N,q,\lambda_{\tilde{r}},a)\inorm{g},\end{split} \end{align} where $K_5$ denotes a arbitrary constant.
\end{enumerate}	

\begin{proof}
	\begin{enumerate}
		\item 
		\begin{align} \label{any3}
		\begin{split}
		\left\| \int_{0}^{\infty} \frac{p}{\tau}~C(t)  \widetilde{B(-\tau)} ~ \pi_ug(\tau) d\tau \right\| \leq  \int_{0}^{\infty} \frac{p}{\tau}~ \mnorm{C(t)}  \mnorm{\widetilde{B(-\tau)}} \inorm{\pi_ug} d \tau \\ 
		\leq \inorm{g} \int_{0}^{\infty} \frac{p}{\tau} ~\mnorm{C(t)}  \mnorm{\widetilde{B(-\tau)}}~ d\tau
		\end{split} 
		\end{align}
		\begin{dmath*}
			\leq \inorm{g} \int_{0}^{\infty} \left( \tilde{K}_1(q,\lambda_{\tilde{r}}) t^{-2p} + \tilde{K}_2(q,\lambda_{\tilde{r}}) t^{-p-pq}\right)~\frac{p}{\tau}~ \left( \sum_{m=0}^{n_{\tilde{r}} -1}\frac{1}{p^{m+1}} M(\lambda_{\tilde{r}},m+1) \tau^{m+1} e^{-\tau \alpha} \right)~ d\tau  
		\end{dmath*}
		\begin{dmath}
			\leq \inorm{g} \left( \tilde{K}_1(q,\lambda_{\tilde{r}}) t^{-2p} + \tilde{K}_2(q,\lambda_{\tilde{r}}) t^{-p-pq}\right) \sum_{m=0}^{n_{\tilde{r}} -1} \frac{1}{p^{m}} M(\lambda_{\tilde{r}},m+1) \left( \int_{0}^{\infty} \tau^m e^{-\tau \alpha} d\tau \right) \leq \inorm{g} \left( \tilde{K}_1(q,\lambda_{\tilde{r}}) t^{-2p} + \tilde{K}_2(q,\lambda_{\tilde{r}}) t^{-p-pq}\right) \sum_{m=0}^{n_{\tilde{r}} -1}\frac{1}{p^{m}} M(\lambda_{\tilde{r}},m+1) \frac{m!}{\alpha^{m+1}}. 
		\end{dmath}
		
Since $t^{-2p}\leq N^{-2p}, t^{-p-pq} \leq N^{-p-pq}$ whenever $t > N$, we have
\begin{align}
\begin{split}
\leq \inorm{g} \bigg[  \tilde{K}_1(q,\lambda_{\tilde{r}}) &\left(  \sum_{m=0}^{n_{\tilde{r}} -1} \frac{1}{p^{m}} M(\lambda_{\tilde{r}},m+1) \frac{m!}{\alpha^{m+1}} \right) N^{-2p}  \\ &+ \tilde{K}_2(q,\lambda_{\tilde{r}}) \left(  \sum_{m=0}^{n_{\tilde{r}} -1} \frac{1}{p^{m}} M(\lambda_{\tilde{r}},m+1) \frac{m!}{\alpha^{m+1}}\right) N^{-p-pq}\bigg].
\end{split}  
\end{align}
		
Denoting the terms in square bracket as $K_3(N,q,\lambda_{\tilde{r}})$ we get
\begin{align}
\left\| \int_{0}^{\infty} C(t) p B(-\tau) \; \pi_ug(\tau) d\tau \right\| \leq K_3(N,q,\lambda_{\tilde{r}}) \inorm{g}.
\end{align}

\item 
\footnotesize{
\begin{align*}	
	&\left\| \int_{0}^{t} (t-\tau)^{p-1} \Gml{p,p}{(t-\tau)^p}{A}~~ \pi_ug(\tau) ~d\tau + \int_{0}^{\infty}\frac{p}{\tau} B(t)\widetilde{B(-\tau)}~~\pi_ug(\tau)~ d\tau\right\| \\
	 \begin{split} &= \bigg|\bigg|\int_{0}^{t-1}(t-\tau)^{p-1} \Gml{p,p}{(t-\tau)^p}{A}~~ \pi_ug(\tau) ~d\tau  \\
	&+ \int_{t-1}^{t}(t-\tau)^{p-1} \Gml{p,p}{(t-\tau)^p}{A}~~ \pi_ug(\tau) ~d\tau  
	+ \int_{0}^{\infty}\frac{p}{\tau} B(t)\widetilde{B(-\tau)}~~\pi_ug(\tau) ~d\tau \bigg|\bigg|
	\end{split}\\
	\begin{split}
&\leq \bigg|\bigg| \int_{0}^{t-1}(t-\tau)^{p-1} \left[ (t- \tau)^{-p} \widetilde{B(t-\tau)} + \widetilde{C(t-\tau)} \right] ~~ \pi_ug(\tau) ~d\tau \\&+ \int_{t-1}^{t}(t-\tau)^{p-1} \Gml{p,p}{(t-\tau)^p}{A}~~ \pi_ug(\tau) ~d\tau + \int_{0}^{\infty}\frac{p}{\tau} B(t)\widetilde{B(-\tau)}~~\pi_ug(\tau) ~d\tau\bigg|\bigg|
\end{split}\\ 
\begin{split}
&\leq \overbrace{\bigg|\bigg| \int_{0}^{t-1} (t-\tau)^{p-1} \widetilde{C(t-\tau)} ~~ \pi_ug(\tau)~ d\tau \bigg|\bigg|}^{\Scale[1.5]{\alpha}} \\&+ \overbrace{\bigg|\bigg| \int_{0}^{t-1} (t-\tau)^{-1}\widetilde{B(t-\tau)} ~~ \pi_ug(\tau) d\tau + \int_{0}^{\infty}\frac{p}{\tau} ~B(t)\widetilde{B(-\tau)}~~\pi_ug(\tau) ~d\tau\bigg|\bigg|}^{\Scale[1.5]{\beta}} \\&+ \underbrace{\bigg|\bigg| \int_{t-1}^{t} (t-\tau)^{p-1} \Gml{p,p}{(t-\tau)^{p}}{A}~~ \pi_ug(\tau) ~d\tau \bigg|\bigg|}_{\Scale[1.5]{\gamma}}.
\end{split}\\
& \leq \Scale[1.5]{\alpha}+\Scale[1.5]{\beta}+\Scale[1.5]{\gamma}  
\end{align*} 
}

The bounds for $\Scale[1.5]{\gamma}$ and $\Scale[1.5]{\alpha}$ are found in \citep[eqn. (70) and eqn. (58)]{deshpande2016local} and remain the same. The proof for the bound on $\Scale[1.5]{\beta}$ is given below.  

Using corrected Lemma 6 for $t > N$,
\begin{align} 
	\nonumber\Scale[1.5]{\beta}&=\bigg|\bigg| \int_{0}^{t-1} (t-\tau)^{-1}~\widetilde{B(t-\tau)} ~~ \pi_ug(\tau) ~d\tau + \int_{0}^{\infty} \frac{p}{\tau} B(t) \widetilde{B(-\tau)} ~~ \pi_ug(\tau) ~d\tau \bigg|\bigg|\\ 
	\nonumber& = \bigg|\bigg| \int_{0}^{t-1} \frac{\widetilde{B(t-\tau)}}{(t-\tau)}~\pi_ug(\tau)~d\tau - \int_{0}^{\infty} \frac{\widetilde{B(t-\tau)}}{(t-\tau)}~\pi_ug(\tau)~d\tau\bigg|\bigg|\\
	&\leq \inorm{g} \left[ \int_{t-1}^{t}\left| \frac{\widetilde{B(t-\tau)}}{(t-\tau)}\right|~ d\tau + \int_{t}^{\infty}\left| \frac{\widetilde{B(t-\tau)}}{(t-\tau)}\right| ~d\tau\right]. \label{err3}
\end{align}
From the asymptotic expansion of $\widetilde{B}$ we have
\begin{align}
\int_{t}^{\infty}\left| \frac{\widetilde{B(t-\tau)}}{(t-\tau)}\right| d\tau	\leq  \int_{t}^{\infty} \sum_{m=0}^{n_{\tilde{r}} -1} \frac{1}{p^{m+1}} M(\lambda_{\tilde{r}},m+1) \frac{(\tau - t)^{m+1}}{(\tau-t)} e^{-(\tau -t) \alpha} ~d\tau,
\end{align} and by integrating terms on right side we get
\begin{align} \label{eq4.7}
\int_{t}^{\infty}\left| \frac{\widetilde{B(t-\tau)}}{(t-\tau)}\right| d\tau	\leq \left( \sum_{m=0}^{n_{\tilde{r}} -1} \frac{1}{p^{m+1}} M(\lambda_{\tilde{r}},m+1) \frac{m!}{\alpha^{m+1}}\right),
\end{align}
and
\begin{align}\label{eq4.8}
\int_{t-1}^{t}\left| \frac{\widetilde{B(t-\tau)}}{(t-\tau)}\right| ~d\tau	\leq \int_{t-1}^{t} \sum_{m=0}^{n_{\tilde{r}}-1} \frac{1}{p^{m+1}}M(\lambda_{\tilde{r}},m+1) (t-\tau)^{m} e^{(t-\tau)\alpha}d \tau,
\end{align}since $t-\tau >0$. Note							 
\begin{align} \label{late1}
\int_{t-1}^{t} (t-\tau)^m e^{(t-\tau)\alpha} d\tau = \int_{0}^{1} u^m e^{u \alpha} du \leq \int_{0}^{1} e^{\alpha} du \leq e^{\alpha}.
\end{align}
In view of eqn. \eqref{late1}, eqn.\eqref{eq4.8} reduces to
\begin{align} \label{eq4.9}
\int_{t-1}^{t}\left| \frac{\widetilde{B(t-\tau)}}{(t-\tau)}\right| ~d\tau	\leq \sum_{m=0}^{n_{\tilde{r}}-1} \frac{1}{p^{m+1}}M(\lambda_{\tilde{r}},m+1) e^{\alpha}.
\end{align}

Substituting eqn. \eqref{eq4.7} and eqn. \eqref{eq4.9} in eqn. \eqref{err3} we get
\begin{align}
\Scale[1.5]{\beta} \leq \left[  \sum_{m=0}^{n_{\tilde{r}}-1} \frac{1}{p^{m+1}}M(\lambda_{\tilde{r}},m+1)\left(  e^{\alpha} + \frac{m!}{\alpha^{m+1}}\right) \right] \inorm{g}.
\end{align}	

Adding \Scale[1.5]{\alpha},\Scale[1.5]{\beta},\Scale[1.5]{\gamma} and renaming the constant as $K_5$ we get
\begin{align}
\begin{split}
\left\| \int_{0}^{t} (t-\tau)^{p-1} \Gml{p,p}{(t-\tau)^p}{A}\; \pi_ug(\tau) d\tau + \int_{0}^{\infty}\frac{p}{\tau} B(t)\widetilde{B(-\tau)}~ \pi_ug(\tau) d\tau\right\|\\ < K_5(N,q,\lambda_{\tilde{r}},a)\inorm{g}.\end{split} \end{align}
Hence the proof.
\end{enumerate}
\end{proof}
Lemma 9 part 1 (cf. \citep[eqn. (74)]{deshpande2016local}) takes the following form while Lemma 9 part 2 and part 3(cf. \citep[eqn. (75) and eqn. (76)]{deshpande2016local}) remain the same.\\
\textbf{Lemma 9 part 1.} Let $g \in C[I, \mathbb{R}^n]$, $t \leq N$ and $K_8(N,a):= C_1 \exp(N a^{1/p} + C_2)$, $C_1,~ C_2$ arbitrary. Further let  $\sigma^{*}: = \int_{0}^{\infty} \frac{p}{\tau}\widetilde{B(-\tau)}~\pi_u g(\tau)~d\tau$. Then
\begin{equation}
\bigg|\bigg|\ml{p}{A} \int_{0}^{\infty} \frac{p}{\tau} \widetilde{B(-\tau)}~\pi_u g(\tau) ~d\tau \bigg|\bigg| \leq \inorm{g} \norm{\sigma^{*}} K_8(N,a).
\end{equation}
Note that from rectified Lemma 8 it is clear that integral in $\sigma^{*}$ exists in $\mathbb{R}^n$, and hence $\sigma^{*}$ is well defined.

The proof follows on similar lines as in the original article \cite{deshpande2016local}.

 Lemma 10 from ref. \cite{deshpande2016local} remains same. For the sake of completeness, we re-state the results below:
 
\textbf{ Lemma 10.}
 	Let $x,y \in C[I,\mathbb{R}^n]$ and $r:=\left\lbrace \inorm{x}, \inorm{y}\right\rbrace $, then
 	\begin{align}
 	\left\| f(x(t)) - f(y(t))\right\|\leq \varepsilon_r \left\| x(t) - y(t) \right\|, \quad (t\geq 0),  
 	\end{align} whenever $x,y \in N_r(0)$.
 	Then for any $\sigma,\widehat{\sigma} \in \mathbb{R}^n$, following inequalities hold:
 	\begin{enumerate}
 		\item \begin{bulletequation}
 			\left\| T_{\sigma}(x) - T_{\widehat{\sigma}}(y)\right\|_{\infty} \leq M_5(N,q,\lambda_{\tilde{r}},a) \varepsilon_r \inorm{x-y} + K_9(N,q,\lambda_{\tilde{r}},a)\norm{\sigma - \widehat{\sigma}}.
 		\end{bulletequation}
 		
 		\item \begin{bulletequation}
 			\left\| T_{\sigma}(x) \right\|_{\infty} \leq M_5(N,q,\lambda_{\tilde{r}},a)\varepsilon_r \inorm{x} + K_9(N,q,\lambda_{\tilde{r}},a) \norm{\sigma}. 
 		\end{bulletequation}
 	\end{enumerate}

 \textbf{Proof of local stable manifold theorem, Step II:} Owing to the changes in the operator $\pi_u$, \citep[eqn. (119) - eqn. (125)]{deshpande2016local} should be replaced by the following.
 
 Consider the unstable projection of $y(t)$ for $t > \tilde{N}$,
 \begin{align} 
\nonumber\begin{split}
 	\norm{\pi_u y(t)} &=  \bigg|\bigg|\ml{p}{A} \int_{0}^{\infty} \frac{p}{\tau}~ \widetilde{B(-
 	\tau)}~~ \pi_u f(y(\tau))~d\tau \\
 &+ \int_{0}^{t} (t-\tau)^{p-1} E_{p,p}((t-\tau)^p A) \pi_u f(y(\tau)) d\tau \bigg|\bigg|
 \end{split}\\
\nonumber\begin{split} 
	 &\leq  \bigg|\bigg| \int_{0}^{\infty} \frac{p}{\tau}~C(t)~ \widetilde{B(-
		\tau)}~~\pi_u f(y(\tau))~d\tau \bigg|\bigg|\\
	&+ \bigg|\bigg|\int_{0}^{t} (t-\tau)^{p-1} E_{p,p}((t-\tau)^p A) \pi_u f(y(\tau)) d\tau \\
	&+ \int_{0}^{\infty} \frac{p}{\tau} B(t)~\widetilde{B(-\tau)}~\pi_uf(y(\tau))~d\tau \bigg|\bigg|
	\end{split}\\ \label{big}
	&\leq \left[ K_3(\tilde{N},q,\lambda_{\tilde{r}}) + K_5(\tilde{N},q,\lambda_{\tilde{r}},a)\right] ~\norm{\pi_u f(y(\tau))},
 \end{align}
 where the last inequality is due to corrected Lemma 8.
 
 Note
 \begin{equation}
 \norm{\pi_u f(y(\tau))} \leq \norm{f(y(\tau))} \leq \varepsilon_{L+\epsilon}(L+\epsilon).
 \end{equation}
 In view of eqn. \eqref{big} 
 \begin{align}
\nonumber\norm{\pi_u y(t)} &\leq \varepsilon_{L+\epsilon}(L+\epsilon)~\left[ K_3(\tilde{N},q,\lambda_{\tilde{r}}) + K_5(\tilde{N},q,\lambda_{\tilde{r}},a)\right] \\
\nonumber&\leq \varepsilon_{L+\epsilon}(L+\epsilon) ~M_5 (\tilde{N},q,\lambda_{\tilde{r}},a)\\
&\leq \xi_{\epsilon} (L + \epsilon).
 \end{align}
 Hence the proof for the step II follows.
 
 The rest of the proof of the theorem remains same.
 \section{Illustrative Examples}
 \label{exmp}
 We discuss the example discussed in ref. \cite{deshpande2016local} below. Due to the corrections proposed in Section 2 of this paper, the local stable manifold will be different, and is presented in the following example.  
 \begin{exmp}
 Consider the following IVP:
 \begin{align}
 	D^p x(t) = A x(t) + f(x(t)),\quad  x(0) = x_0 = \begin{pmatrix}
 		\sigma_1 \\ \sigma_2 \\ \sigma_3
 	\end{pmatrix},
 \end{align}
 where $x(t) = \begin{pmatrix}
 x_1(t) \\ x_2(t) \\ x_3 (t) \end{pmatrix}$, $f(x) = \begin{pmatrix}
 0 \\ (x_1(t))^2 \\ 3 (x_1(t))^2
 \end{pmatrix}$ and $
 A = \begin{pmatrix}
 -1 & 0 & 0 \\ 
 0 & 2 & 1 \\ 
 0 & 0 & 2
 \end{pmatrix}.$\\ 
 By solving linear system $D^p x = A x$, we find
 \begin{align}
 	E^s = \{x_0 \in \mathbb{R}^n \slash \sigma_2 = \sigma_3 = 0 \}.
 \end{align} 
 The stable and unstable projections are defined as    
 \begin{align}
 	\pi_u x = \begin{pmatrix}
 		0 \\ x_2 \\ x_3
 	\end{pmatrix}, \quad \pi_s x = \begin{pmatrix}
 	x_1 \\ 0 \\0
 \end{pmatrix}, \text{ where  }x = (x_1,x_2,x_3)^T \in \mathbb{R}^3\end{align}
It may be noted that
 \begin{align} \label{ex2}
 	x_1(t) = E_p(-t^p) \sigma_1.
 \end{align}
 By using unstable projection of $T_{\sigma}$ where $\sigma = (\sigma_1,0,0)$ and using the fact that $x(t)$ is a fixed point of $T_{\sigma}$, we get
 \begin{align} \label{ex3}
 	\begin{split}
 		\pi_u x(t) = \ml{p}{A} &\int_{0}^{\infty} \frac{p}{\tau} \widetilde{B(-\tau)}~ \pi_u f(x(\tau))~ d\tau \\&+ \int_{0}^{t} (t-\tau)^{p-1} E_{p,p}((t-\tau)^p A)~ \pi_u f(x(\tau))~ d\tau.
 	\end{split}
 \end{align}
 
 Note \begin{align}\label{ex11}
 	\pi_uf(x(\tau)) = \begin{pmatrix}
 		0 \\ x_1^2(\tau) \\ 3 x_1^2 (\tau)
 	\end{pmatrix} = \begin{pmatrix}
 	0 \\ E_{p}^{2}(-\tau^p)\sigma_1^2 \\ 3 E_{p}^{2}(-\tau^p)\sigma_1^2 
 \end{pmatrix}
 \end{align} 
 and 
 \begin{align} \label{ex12}
 	\frac{p}{\tau}~ \widetilde{B(-\tau)} &= \begin{pmatrix}
 		0 & 0 & 0 \\ 0 & \frac{p}{\tau}\Psi_{0}(-\tau,2) & \frac{p}{\tau}\Psi_1(-\tau,2) \\ 0 & 0 & \frac{p}{\tau}\Psi_{0}(-\tau,2)
 	\end{pmatrix}.\\
 	\nonumber
 	&= \begin{pmatrix}
 	0 & 0 & 0 \\ 0 & -2^{\frac{1}{p}-1} e^{-2^{\frac{1}{p}}\tau} & \frac{1}{p} 2^{\frac{1}{p}-2} e^{-2^{\frac{1}{p}}\tau} ~(2^{\frac{1}{p}}\tau -1 + p) \\
 	0 & 0 & -2^{\frac{1}{p}-1} e^{-2^{\frac{1}{p}}\tau}
 	\end{pmatrix}
 \end{align}
 In view of the values of $\pi_u f(x(\tau))$ and $\frac{p}{\tau} \widetilde{B(-\tau)}$ given in eqn. \eqref{ex11}, eqn. \eqref{ex12} respectively,  eqn. \eqref{ex3} implies
 \begin{dmath} \label{ex4}
 	x_3(t) = -3 E_p(2 t^p) ~ \mathbf{l} ~ \sigma_1 ^2 ~2^{\frac{1}{p}-1} + 3 \sigma_1^2 \int_{0}^{t} (t-\tau)^{p-1} E_{p,p}(2(t-\tau)^p) E^2_p(-\tau^p)~ d\tau,
 \end{dmath}
 and
 \begin{align}\label{ex5}
 \begin{split}
 	x_2(t) &= - E_P(2 t^p) ~\mathbf{l}~\sigma_1^2~2^{\frac{1}{p}-1} + \frac{3}{p}~ E_p(2 t^p)~ \mathbf{m}~ \sigma_1^2~2^{\frac{1}{p}-2} \\&- \pder{}{\lambda}\ml{p}{\lambda}\bigg\vert_{\lambda =2} (3 ~\mathbf{l} ~\sigma_1^2) ~ 2^{\frac{1}{p}-1}  \\&+\sigma_1^2 \int_{0}^{t} (t-\tau)^{p-1} E_{p,p}(2(t-\tau)^p) E_p^2(-\tau^p)~ d\tau \\&+ 3\sigma_1^2 \int_{0}^{t} (t-\tau)^{p-1} \left( \pder{}{\lambda} E_{p,p}((t-\tau)^p \lambda)\right) \bigg\vert_{\lambda=2} E_p^2(-\tau^p) d\tau,
 	\end{split}
 \end{align}
 where \begin{align*}
 	l &= \int_{0}^{\infty} e^{-\tau 2^{1/p}} E_p ^2(-\tau^p) d\tau,\\
 	m &= \int_{0}^{\infty}  e^{-\tau 2^{1/p}}~(\tau 2^{\frac{1}{p}} -1 + p)~ E_p^2(-\tau^p)~d\tau. 
 \end{align*} 
 For sufficiently small neighborhood of origin and $t=0$, eqn. \eqref{ex4} and eqn. \eqref{ex5} yield the map
 \begin{align}
 	\sigma_3 &= x_3(0) = -3 \;l\;\sigma_1^2~2^{\frac{1}{p}-1},\\
 	\sigma_2 &=x_2(0) = -l~ \sigma_1^2~2^{\frac{1}{p}-1} + \frac{3}{p}~m \sigma_1^2~ 2^{\frac{1}{p}-2},
 \end{align}
which defines the required local stable manifold over $E^s$.   
\end{exmp}

The following example was discussed in \cite{cong2016stable}. We point out that this example is in agreement with the main result of the present article.
\begin{exmp}
	 Consider the IVP:
	 \begin{align}
	 D^p x(t) = A x(t) + f(x(t)),\quad  x(0) = x_0 = \begin{pmatrix}
	 \sigma_1 \\ \sigma_2
	 \end{pmatrix},
	 \end{align}
	 where $x(t) = \begin{pmatrix}
	 x_1(t) \\ x_2(t)  \end{pmatrix}$, $f(x) = \begin{pmatrix}
	 (x_{1}(t))^2 \\ (x_1(t))^2 + (x_2(t))^2 
	 \end{pmatrix}$ and $
	 A = \begin{pmatrix}
	 -2 & 0 \\
	 0 & 2
	 \end{pmatrix}.$\\ 
	 Note that $E^s = \left\lbrace x_0 \in \mathbb{R}^2 / \sigma_2 = 0\right\rbrace $.
	 The stable and unstable projections are defined as    
	 \begin{align}
	 \pi_u x = \begin{pmatrix}
	 0 \\ x_2 
	 \end{pmatrix}, \quad \pi_s x = \begin{pmatrix}
	 x_1 \\ 0 
	 \end{pmatrix}, \text{ where  }x = (x_1,x_2)^T \in \mathbb{R}^2.\end{align}
	 By using unstable projection of $T_{\sigma}$ where $\sigma = (\sigma_1,0)$ and using the fact that $x(t)$ is a fixed point of $T_{\sigma}$, we get
	 \begin{align} \label{ex23}
	 \begin{split}
	 \pi_u x(t) = \ml{p}{A} &\int_{0}^{\infty} \frac{p}{\tau} \widetilde{B(-\tau)}~ \pi_u f(x(\tau))~ d\tau \\&+ \int_{0}^{t} (t-\tau)^{p-1} E_{p,p}((t-\tau)^p A)~ \pi_u f(x(\tau))~ d\tau.
	 \end{split}
	 \end{align}
	 
	 Let $f(x(\tau)) = \phi(\tau)$, we have
	  \begin{align}\label{ex211}
	 \pi_uf(x(\tau)) = \begin{pmatrix}
	 0 \\ (\phi_1(\tau))^2 + (\phi_2(\tau))^2
	 \end{pmatrix} 
	 \end{align} 
	 and 
	 \begin{align} \label{ex212}
	 \frac{p}{\tau}~ \widetilde{B(-\tau)} &= \begin{pmatrix}
	 0 & 0  \\ 0 & \frac{p}{\tau}\Psi_{0}(-\tau,2) 
	 \end{pmatrix}= \begin{pmatrix}
	 0 & 0  \\ 0 & -2^{\frac{1}{p}-1} e^{-2^{\frac{1}{p}}\tau} 
	 \end{pmatrix}.
	 \end{align}
	 In view of the values of $\pi_u f(x(\tau))$ and $pB(-\tau)$ given in eqn. \eqref{ex211}, eqn. \eqref{ex212}, second component of eqn. \eqref{ex23} gives
	 \begin{dmath} \label{ex24}
	 	x_2(t) = -E_p(2 t^p)~ \int_{0}^{\infty}  ~2^{\frac{1}{p}-1}~e^{-\tau 2^{\frac{1}{p}}}~(\phi_1^2 + \phi_2^2)~d\tau +  \int_{0}^{t} (t-\tau)^{p-1} E_{p,p}(2(t-\tau)^p) ~ (\phi_1^2 + \phi_2^2)~d\tau.
	 \end{dmath}
	 We claim that $\norm{\pi_u x} \rightarrow 0$ as $t \rightarrow \infty$ provided $\norm{\phi(t)} \rightarrow 0$ as $t \rightarrow \infty$.
	 
	 Consider 
	  \begin{dmath} \label{bigg} 
	  	\norm{\pi_u x} = \norm{-E_p(2 t^p)~ \int_{0}^{\infty}  ~2^{\frac{1}{p}-1}~e^{-\tau 2^{\frac{1}{p}}}~(\phi_1^2 + \phi_2^2)~d\tau +  \int_{0}^{t} (t-\tau)^{p-1} E_{p,p}(2(t-\tau)^p) ~ (\phi_1^2 + \phi_2^2)~d\tau} \leq 
	  	 \norm{ \int_{0}^{\infty}  ~-2^{\frac{1}{p}-1}~e^{-\tau 2^{\frac{1}{p}}} ~C(t)~(\phi_1^2 + \phi_2^2)~d\tau} +  \norm{\int_{0}^{t} (t-\tau)^{p-1} E_{p,p}(2(t-\tau)^p) ~ (\phi_1^2 + \phi_2^2)~d\tau + \int_{0}^{\infty} B(t) (-2^{\frac{1}{p}-1}~e^{-\tau 2^{\frac{1}{p}}})~(\phi_1^2 + \phi_2^2)~d\tau}
	  \end{dmath}
	  From corrected Lemma 8 part 1 and 2, for sufficiently large $t$ eqn. \eqref{bigg} becomes
	  \begin{equation} \label{end}
		\norm{\pi_u x} \leq K~ (\phi_1^2 + \phi_2^2),
	  \end{equation} where $K$ is an arbitrary constant.
	  Since $\norm{\phi(t)} \rightarrow 0$ is given, eqn. \eqref{end} proves the claim. 
\end{exmp}

\begin{exmp}
	The fractional ordered Liu system is defined as
	\begin{align}
	\begin{split} \label{liu}
	D^{\alpha}x_1 &= -a x_1 - e x_2^2 \\
	D^{\alpha}x_2 &= b x_2 - k x_1 x_3 \\
	D^{\alpha}x_3 &= -c x_3 + m x_1 x_2, 
	\end{split}
	\end{align}
	where parameter values are taken as $a=1, e=0, b=2.5, k=4, c=5, m=0$. Let $x(t) = \begin{pmatrix}
	x_1(t) \\ x_2(t) \\ x_3(t) 
	\end{pmatrix}, f(x) =  \begin{pmatrix}
	0 \\ -4 x_1(t) x_3(t) \\ 0 
	\end{pmatrix}$ and $A = \begin{pmatrix}
	-1 & 0 & 0 \\ 0 & 2.5 & 0 \\ 0 & 0 & -5 
	\end{pmatrix}$ where matrix $A$ denotes the Jacobian matrix of the system \eqref{liu} around an equilibrium point $(0,0,0)^T$. Let $ x(0) =x_0= \begin{pmatrix}
	\sigma_1(t) \\ \sigma_2(t) \\ \sigma_3(t) 
	\end{pmatrix}$ be the given initial condition. In this case system \eqref{liu} can be re-written in the form
	\begin{equation}
	D^{\alpha} x = A x + f(x), ~x(0) = x_0.
	\end{equation} 
	
	Solving the linear part of the system \eqref{liu} we get stable and unstable subspaces as
	\begin{align}
	E^s = \left\lbrace x_0 \in \mathbb{R}^3 \bigg/ \sigma_2 = 0\right\rbrace, ~~ E^u = \left\lbrace x_0 \in \mathbb{R}^3 \bigg/ \sigma_1 = \sigma_3 = 0 \right\rbrace.
	\end{align}
	Thus the stable and unstable projection maps are given as 
	\begin{align}
	\pi_u x = (0, x_2, 0)^T, ~~	\pi_s x = (x_1, 0, x_3)^T.
	\end{align}
	
	Solving system \eqref{liu} we get
	\begin{equation}
	x_1(t) = \mll{p}{-} \sigma_1, ~\text{and}~ x_3(t) = \mll{p}{-5} \sigma_3. 
	\end{equation}
	Using unstable projection of operator in \eqref{define} and fact that $x$ is the fixed point we get
	\begin{align}
	\begin{split} \label{liu1}
	\pi_u x(t) = \ml{p}{A} \int_{0}^{\infty} \frac{\alpha}{\tau} &\widetilde{B}(-\tau)~ \pi_u f(x(\tau)) ~d\tau \\&+ \int_{0}^{t} (t-\tau)^{p -1} E_{p,p}((t-\tau)^{p} A)~\pi_u f(x(\tau)) ~d\tau.
	\end{split}
	\end{align}
	 Now $\pi_u f(x(\tau))  = (0, -4 \mll{p}{-} \mll{p}{-5} \sigma_1 \sigma_3,0)^T$.
	 Further
	 \begin{align}
	 \frac{p}{\tau} \widetilde{B}(-\tau) = \begin{pmatrix}
	 0 & 0 & 0 \\ 0 & \frac{p}{\tau} \widetilde{\Psi_0}(-\tau,5/2) & 0 \\ 0 & 0 & 0
	 \end{pmatrix} = \begin{pmatrix}
	 0 & 0 & 0 \\ 0 & - (\frac{5}{2})^{\frac{1}{p}-1} ~e^{-\tau (5/2)^{1/p}}& 0 \\ 0 & 0 & 0
	 \end{pmatrix}.
	 \end{align}
	 In light of this, equation \eqref{liu1} gives
	 \begin{align}
	 \begin{split}\label{liu2}
	 x_2(t) = \mll{p}{\frac{5}{2}} ~\int_{0}^{\infty} &4~ (\frac{5}{2})^{\frac{1}{p}-1}~e^{-\tau (5/2)^{1/p}} ~E_p(-\tau^p)~E_p(-5 \tau^p) \sigma_1 \sigma_3~ d\tau \\
	 &-4 \sigma_1 \sigma_3~ \int_{0}^{t}  (t-\tau)^{p -1} E_{p,p}(\frac{5}{2}(t-\tau)^{p})~E_p(-\tau^p)~E_p(-5 \tau^p) d\tau.
	 \end{split}
	 \end{align}
	 	 
	 For sufficiently small neighborhood of origin and $t=0$, eqn. \eqref{liu2} yields the map
	 \begin{equation}\label{liu3}
	 \sigma_2 = x_2(0) = l ~\sigma_1 ~\sigma_3, \text{ where}
	 \end{equation}
	  \begin{equation}
	 l = 4~ (\frac{5}{2})^{\frac{1}{p}-1} ~\int_{0}^{\infty} e^{-\tau (5/2)^{1/p}} ~E_p(-\tau^p)~E_p(-5 \tau^p) ~ d\tau.
	 \end{equation}
	 For $p=0.5$, the local stable manifold of the \eqref{liu} around origin is plotted in figure (\ref{figliu}).
	\begin{figure}[H]
		\centering
		\includegraphics{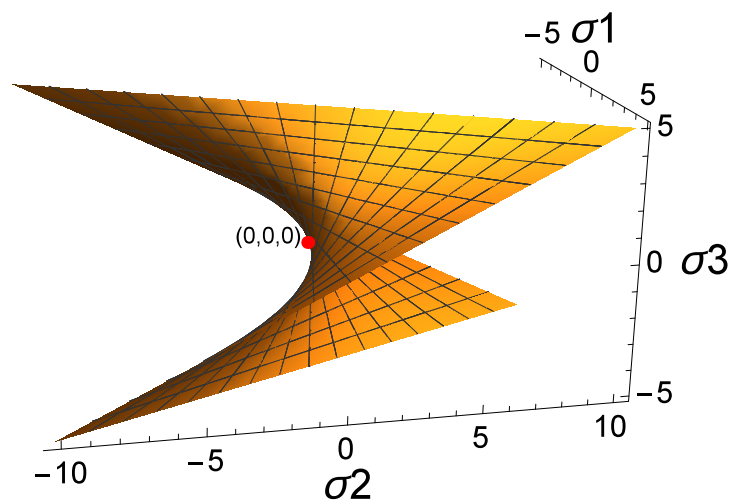}
		\caption{ local stable manifold around equilibrium point $(0,0,0)$ of fractional Liu system for $p=0.5$.}
			\label{figliu}
	\end{figure}
\end{exmp} 
\section{Conclusions}
	Equation in Lemma 4 part 2 of the ref. \cite{deshpande2016local} has been corrected. Further local stable manifold theorem has been established following the approach given in \cite{deshpande2016local}. The example given in \cite{cong2016stable} has been discussed and it is shown that $\pi_u x \rightarrow 0$ in view of the corrections presented here. 
 	\bibliographystyle{ieeetr}

\end{document}